\documentclass[11pt,a4paper]{article}
\usepackage{amsmath}
\usepackage{amssymb}
\usepackage{theorem}
\usepackage{xspace}

\def\pth#1{\left(#1\right)}
\def\acc#1{\left\{#1\right\}}
\def\cro#1{\left[#1\right]}

\def\eR{I\!\!R}
\def\eE{I\!\!E}
\def\eP{I\!\!P}
\def\e1{1\!\!1}
\def\eN{I\!\!N}
\def\el{{\lambda}'}
\def\ef{\overset{.}{f}}

\newtheorem{theorem}{Theorem}[section]

\newtheorem{lemma}{Lemma}[section]

\begin{document}

\title{The  M-estimation in a multi-phase random nonlinear model }

\author{Gabriela CIUPERCA\footnote{Universit\'e de Lyon, Universit\'e Lyon 1, 
UMR 5208, blvd du 11 novembre 1918,
F - 69622 Villeurbanne Cedex, France, E-mail:Gabriela.Ciuperca@univ-lyon1.fr}}


\maketitle
\begin{abstract}
This paper considers M-estimation of a nonlinear regression model with multiple change-points occuring at unknown times. The multi-phase random design  regression model, discontinuous in each change-point, have an arbitrary  error $\varepsilon$. In the case when the number of jumps is known, the M-estimator of  locations of breaks and of regression parameters are studied. These estimators are consistent and the distribution of  the regression parameter  estimators is  Gaussian. The estimator of each   change-point converges, with the rate $n^{-1}$, to the smallest minimizer  of the  independent compound Poisson processes. The results are valid for a large class of error distributions.

\end{abstract}

\section{Introduction}
Change-points are intrinsic features of signals that appear in economics, medicine and physical science. The statistics literature contains a vast amount of works on issues related to the estimation of the change-point for a parametric regression, most of it specifically designed for the case of a single break. The more used estimators are the maximum likelihood estimators, the least squares estimators or a wider class, the M-estimators. Statistical inference for a parametrical model is influenced by the continuity or by discontinuity of the regression function at the change-points, but also by the determinist character or not of the explicative variable.  We give a non-exhaustive  list with the recent papers. The area of research is so active that it is nearly impossible to list all the recent papers written.\\
For the least squares (LS) estimators we refer to Feder  (1975a, 1975b) for continuous two-lines models, Lai et al. (1979), Yao and Au (1989) for a step function, Liu et al. (1997), Bai and Perron  (1998)   for multiple structural changes in a linear model.\\
For the maximum likelihood (ML) estimator, when the design is determinist, Bhattacharya (1994) discusses his limiting behaviour for a discontinuous linear model. Gill (2004), Gill and Baron (2004) consider a  model where the canonical parameter of an exponential family gradually begins to drift from its initial value at an unknown change-point. For a random design we refer to  Koul and Qian (2002) for two lines model, Ciuperca (2004) for a single jump in a nonlinear model, Ciuperca and Dapzol  (2008) for multiple change-points in linear and nonlinear model.  If the model variance depends of the mean,  the quasi-likelihood estimator can be considered. Braun et al. (2000) consider that the mean is constant between two change-points. Chiou and Muller (2004)  propose a semi-parametric estimator in a  generalized linear model with determinist design.\\
In the general case of  M-estimators,  Rukhin and Vajda (1997) consider  the change-point estimation problem as a nonlinear regression problem, the model being continuous, with a single change-point and fixed design. Koul et al. (2003) study the M-estimators in two-phase linear regression with random design.\\
The present paper makes several contributions to the existing literature. The considered design is  random, the regression function is nonlineary within the framework of a multi-regime and not lastly, a general method of estimation. We study the properties of the M-estimator in a multi-phase discontinuous nonlinear  random regression model with a general error distribution. The class of the M-estimators was introduced by Huber (1964) and its principal properties are exposed in Huber (1981). We generalize among others, the results for the two-phase random linear  model of Koul et al. (2003) obtained  by M-estimation, the results obtained  by the ML estimation of  Ciuperca and Dapzol (2008) for a multiphase random nonlinear model and of Bai and Perron (1998) obtained  by  LS estimation in a multiple nonrandom linear regression.  An important point of the proofs for the linear case is the relation between the regression function and its derivatives with respect to regression parameters. Thus we have to modify the approach for the non linear regression. Also, in the case of a single change-point, each of two regimes has one fixed boundary. For multiple breaks, each middle regime has boundaries completely unknown. \\
The paper is organized as follows. We give necessary notations and definitions in Section \ref{Section 2}. In Section \ref{Section 3} we establish the estimators consistency  and the convergence rate. Weak convergence results are also obtained: the asymptotic distribution of the regression parameters M-estimator is Gaussian.   We also  prove that $n(\hat \theta_{2n}-\theta^0_2)$ converges weakly to the smallest minimizer vector of the independent  compound Poisson processes, where $\hat \theta_{2n}$ is the change-point estimator.  Auxiliary results are given  in Appendix.
\section{Notations and model}
\label{Section 2}
Consider the step-function with $K$ ($K \geq 1$) fixed change-points, for $x \in \eR$:
$$
f_\theta(x)=h_{\alpha_0}(x)\e1_{x \leq \tau_1}+h_{\alpha_1}(x)\e1_{\tau_1<x \leq \tau_2}+.....+h_{\alpha_K}(x)\e1_{\tau_K<x }
$$
where $\theta_1=(\alpha_0,\alpha_1,....,\alpha_K) $ are  the nonlinear  regression parameters  and $\theta_2=(\tau_1,...,\tau_K)$, $\tau_1 < \tau_2 < ... < \tau_K$ are the change-points. 
 For all $k=0,1,...,K$, we have the parameter $\alpha_k$ belongs to some compact $\Gamma \subseteq \eR^{d }$. We consider that the vector $\theta_2 \in \eR^K$ and we set $\theta=\pth{\theta_1,\theta_2} \in  \Omega=\Gamma^{K+1}\times  \eR^K$.\\
Consider the random  design model:
$$
Y_i=f_\theta(X_i)+\varepsilon_i, \qquad i=1,...,n
$$
where $(\varepsilon_i,X_i)$ is a sequence of continuous independent random variables with the same joint distribution as  $(\varepsilon,X)$. The parameter $\theta_1$ and the change-points (or break points) are unknown. The purpose is to estimate $\theta=(\theta_1,\theta_2)$ when $n$ observations of $(Y,X)$ are available.\\
We denote the true value of a parameter with a 0 superscript. In particular, $\theta^0_1=(\alpha^0_0,\alpha^0_1,...,\alpha^0_K)$ and $\theta^0_2=(\tau^0_1,...,\tau^0_K)$ are used to denote, respectively, the true values of the regression parameters and the  true change-points.  Let be also $\theta^0=(\theta^0_1,\theta^0_2)$.  
 We suppose that $\theta^0_1$ is an inner point of the set $\Gamma^{K+1}$.\\
The random variables $X$ and $\varepsilon$ satisfy the following assumptions:\\
{\bf (A1)} $X$ has a positive absolutely continuous Lebesgue density $\varphi$ on $\eR$. Moreover,   $\eE(X^2)< \infty$;\\
{\bf (A2)} $\varepsilon$ has a density absolutely continuous and positive everywhere on $\eR$. Moreover, $\eE(\varepsilon)=0$, $\eE(\varepsilon^2)< \infty$;\\
{\bf (A3)} the random variables $X_i$ and $\varepsilon_i$ are independent.\\
In  the case of linear model with a single change-point: $h_\alpha(x)=a+bx$, $\alpha=(a,b)$  and $K=1$,  assumptions (A1)-(A3)  on $X$ and $\varepsilon$ are also  considered by   Koul et al (2003). \\

\noindent The nonlinear 
 function $h_\alpha$ satisfies the conditions:\\
{\bf (B1)} for all $x \in \eR$, $h_\alpha(x)$ is  three times differentiable with respect to $\alpha$; \\
{\bf (B2)} for  all $x \in \eR$, $\| \partial h_{\alpha_0}(x)/\partial \alpha \| \neq 0 $;\\
{\bf (B3)}  the derivatives $ \partial^3 h_\alpha(x)/\partial \alpha^3$, exist for $ x \in \eR$ and there exist functions $ F_0,F_1,F_2 \in L^2(\varphi)$ such that:
\begin{equation}
\label{eth}
\sup_{\alpha \in \Gamma}|h_\alpha(x)| \leq F_0(x), \qquad \sup_{\alpha \in \Gamma}\|\partial^jh_\alpha(x)/\partial \alpha^j \| \leq F_j(x), \;\; j=1,2
\end{equation}
Obviously, in the case $h_\alpha(x)=a+bx$, the assumptions (B1), (B2) are verified  and (B3) is transformed in (A1). If $h_\alpha(x)$ is a polynomial with degree $p$,  assumption (B3) can be replaced by $\eE(X^{p+1}) < \infty$. \\
Assumption (B2) is necessary for  obtaining the convergence rate of regression parameters estimator.\\

\noindent Let us consider the functions: $d_{(\alpha_k,\alpha_{j})}(x):=h_{\alpha_k}(x)-h_{\alpha_{j}}(x)$, $x \in \eR$, $k,j \in \{1,...K\}$ and the jump at the true break point: $d_k^0:=d_{(\alpha_k^0,\alpha^0_{k-1})}(\tau^0_k)$. We make the identifiability assumption that the jump at each $\tau^0_k$ is non-zero:
\begin{equation}
\label{cidn}
d_{(\alpha_k,\alpha_{k+1})}(\tau^0_k) \neq 0 , \qquad \forall \alpha_{k},\alpha_{k+1} \in \Gamma, \alpha_k \neq \alpha_{k+1} 
\end{equation}
a condition which implies that the function $f_\theta$ is not continuous in the true break points for all parameters in $\Gamma$. For $\theta^*=(\theta^*_1,\theta^*_2)$ and  $\theta=(\theta_1,\theta_2)$, let us denote by $\delta_{(\theta,\theta^*)}(x):=f_\theta(x)-f_{\theta^*}(x)$ the difference between two models. Note also: $\overset{.}f_\theta(x)=\partial f_\theta(x)/\partial \theta_1$.\\

 In the following, we denote by $C$ a generic positive finite constant not depending on $n$. \\
For a vector, let us denote by $\|.\|$ the Euclidean norm and for a matrix $A=(a_{ij})$, $\|A\|={\underset{i,j}{\sum}}|a_{ij}|$. For a vector $v=(v_1,...,v_K)$ we make the convention that $|v|=(|v_1|,...,|v_K|)$.\\

The most important method of constructing statistical estimators is to choose the estimator to maximize or minimize a certain criterion function. The such estimators are called the M-estimators. The maximum likelihood (ML), least squares, least absolute deviation estimators are particular cases. For a function $\rho: \eR \rightarrow \eR_+$, let the M-process be:
$$
M_n(\theta)=\sum^n_{i=1}\rho\pth{Y_i-f_\theta(X_i)}
$$
The following assumptions are  considered for the function $\rho$:\\
{\bf (C1)}  $\rho$ is convex on $\eR$ with right-continuous non-decreasing almost everywhere derivative $\psi$ satisfying $\eE_\varepsilon[\psi^2(\varepsilon+y)]<\infty$, $\forall y \in \eR$. The function $\lambda(y):=\eE_\varepsilon[\psi(\varepsilon+y)]$, $y \in \eR$, is strictly  increasing on $\eR$ and $\lambda$ is continuous at 0 with $\lambda(0)=0$.\\
{\bf (C2)} for all $c \in \eR$, where $\bar \Omega$ is the closure of $\Omega$.\\
{\bf (C3)} the function $y \rightarrow \eE[|\psi(\varepsilon+c+y)-\psi(\varepsilon) |]$ is continuous at 0, $\forall c \in \eR$.\\
{\bf (C4)} the function $\lambda$ is differentiable in a neighborhood of 0, with  derivative $\el$ satisfying $\el(0) \neq 0$, and $\lim_{a \rightarrow 0} a^{-1} \int^a_0 |\el(s)-\el(0)|ds=0$.\\
{\bf (C5)} the random variables $\rho(\varepsilon \pm d_k^0)-\rho(\varepsilon)$, $\forall k=1,...,K$, are continuous.\\

\noindent Assumptions (C1), (C2) are necessary for obtaining the consistency of the estimators, while (C1)-(C5) are used for obtaining the rate of convergence and the asymptotic distribution.\\
Notice that forthe  two-phase linear regression function: $f_\theta(x)=(a_0+b_0x) \e1_{x \leq \tau}+(a_1+b_1x) \e1_{x > \tau}$,  Koul et al. (2003) consider the same assumptions (C1)-(C5). Obviously, (C2) becomes: $\eE_{(\varepsilon,X)}[\psi^2+c_1+c_2|X|] < \infty$, $\forall c_1,c_2 \in \eR$ and  (C5) becomes: $\rho(\varepsilon \pm d) - \rho(\varepsilon)$ continuous, with $d=(a^0_1-a^0_0)+\tau^0(b^0_1-b^0_0)$, $(a^0_0,b^0_0,a_1^0,b_1^0,\tau^0)$ the true value of  $(a_0,b_0,a_1,b_1,\tau)$.\\
 
\noindent For each $\eta>0$, denote the $\eta$-neighborhood of $\theta \in  \Omega$ by:$$
\Omega_\eta(\theta):=\acc{\theta^* =\pth{\theta^*_1,\theta^*_2} \in \Omega \;\; / \;\; \|\theta^*_1 - \theta_1 \| \leq \eta , \|\theta^*_2 -\theta_2\| \leq \eta}
$$
The M-estimator is defined by:
$$
\hat \theta_n:=\pth{\hat \theta_{1n},\hat \theta_{2n}}=\arg \min_{\theta \in \bar \Omega} M_n(\theta) \quad \textrm{ a.s.} 
$$
where $\bar \Theta$ is the close of $\Theta$. Let $\bar \eR=\eR \cup \{-\infty, \infty\}$. The set $\bar \eR$ is compact under the metric $m(x,y)=|arc \tan x - arc\tan y|$, $x,y \in \eR$. For constructing the M-estimator, first we search the regression parameters estimator and then we localize the change-points. First, for a given  $\theta_2 \in \eR^K$, we set:
$$
\tilde \theta_{1n}(\theta_2):=\arg \min_{\theta_1 \in \Gamma^{K+1}} M_n(\theta_1,\theta_2)
$$
Since the number $K$ of the change-points is fixed, the estimator $\tilde \theta_{1n}(\theta_2)$ is constant in $\theta_2$ over any interval of two consecutive ordered $X_i$'s. The M-process $M_n(\tilde \theta_{1n}(\theta_2),\theta_2)$ has only a finite number of possible values with  change-points located at the ordered  $X_{i}$'s. Second, we find the minimizer $\tilde \theta_{2n}$ of $ M_n(\tilde \theta_{1n}(\theta_2),\theta_2)$ with respect to $\theta_2$ over the sample percentile $\{X_i, i=1,...,n\}$. This minimizer may be taken as the left end point of the interval over which  it is obtained.  Then  $\tilde \theta_{2n}=\hat \theta_{2n}$ and the M-estimator is: $\hat \theta_n=(\tilde \theta_{1n}(\tilde \theta_{2n}),\tilde \theta_{2n})$. \\

{\bf Remark.} The considered  model and the estimator are very general. The class of M-estimators includes the least squares ($\rho(x)=x^2$), maximum likelihood ($\rho(.)=\log \varphi_\varepsilon(.)$, with $\varphi_\varepsilon$ the density of $\varepsilon$) and least absolute deviations estimators ($\rho(x)=|x|$).  Examples of distributions satisfying these conditions include Normal for  $X$, double exponential or Normal for the errors  $\varepsilon$ if  $\rho(x)=|x|^a$, $a \in \{1/2,2\}$. \\
For the ML  estimator in a multi-phase nonlinear random model, the conditions imposed on the random variables are (A1), (A2), (A3) and the density $\varphi_\varepsilon$ of $\varepsilon$  satisfies: $u(x)=\varphi'_\varepsilon(x)/\varphi_\varepsilon(x)$ the score function is Holder, is differentiable and $u'$ is Holder also (see Ciuperca and Dapzol(2008)). The function $h_\alpha$ satisfy the condition (B2) and:
$$
\sup_{\theta,\theta^* \in \Omega} \eE_{(\varepsilon,X)}[u^2(\varepsilon+f_\theta(X)-f_{\theta^*}(X))] < \infty
$$ 
\section{Asymptotic properties }
\label{Section 3}
In this section we focus on study of the asymptotic properties of estimator. First, we study the convergence of the M-estimator and we find the rate of convergence. 
\subsection{Consistency and rate of convergence}

For each change-point $\tau^0_k$, since the density of $X$ is absolutely continuous in $\eR$, we have:
\begin{equation}
\label{n2}
n^{-1} \sum^n_{i=1} \e1_{|X_i-\tau^0_k| \leq B/n}=O_{\eP}(n^{-1})
\end{equation}
It is interesting to mention that in a identifiable regular model for a density with jumps, the ML estimator is of order $n^{-1}$ (see Ibragimov and Has'minskii (1981)). For the multi-phase problem, we obtain that the M-estimator of the change-point  has the same order of convergence. Always in a regular model, van der Vaart and Wellner (1996) obtain the rate of convergence of the M-estimator.\\

\noindent The next theorem establishes the strong consistency of the M-estimator and shows that the rate of the convergence of $\hat \theta_{2n}$ to $\theta^0_2$ is $n^{-1}$ and $n^{-1/2}$ of $\hat \theta_{1n}$ to $\theta^0_1$. The theorem includes the results derived by Koul et al. (2003) when $h_\alpha$ is linear for the M-estimator and by Ciuperca and Dapzol (2008) when $h_\alpha$ is nonlinear for the ML estimator. Remark that for the ML estimator in a nonlinear random model, the discontinuity of the function $f_\theta(x)$ in the change-points is not necessary to show the consistency of the estimators.\\
\noindent In order to simplify the study of the rate of convergence,  three processes defined as the differences between two M-processes are considered. The first one is  the difference between a M-process calculated in a some point $\theta$ and a M-process at the true point $\theta^0$:
\begin{equation}
 \label{dntheta}
D_n(\theta_1,\theta_2):=M_n(\theta_1,\theta_2)-M_n(\theta^0_1,\theta^0_2)
\end{equation}
For the second one, the regression parameters vary around $\theta^0_1$,   for $w_1 \in \Gamma^{K+1}$:
$$
D^{(1)}_n(w_1):=M_n\pth{\theta^0_1+n^{-1/2}w_1,\theta^0_2}-M_n(\theta^0_1,\theta^0_2)
$$
the coefficient of $w_1$ being the rate of convergence of the estimator $\hat \theta_{1n}$ and finally we make vary the change-points:
$$
D^{(2)}_n(\theta_1,\theta_2):=M_n(\theta_1,\theta_2)-M_n(\theta_1,\theta^0_2)
$$
The relation between these processes is given by the following decomposition:
\begin{equation}
\label{eq3.13}
D_n(\theta_1,\theta_2)=D^{(1)}_n\pth{n^{1/2}(\theta_1-\theta^0_1)}+D^{(2)}_n(\theta_1,\theta_2)
\end{equation}
\begin{theorem}
\label{theorem 1}{\it (i)} Under  assumptions (\ref{cidn}), (A1), (A3), (B1), (B3), (C1) and (C2) we have:  $
\hat \theta_n \overset{a.s.} {\underset{n \rightarrow \infty}{\longrightarrow}} \theta^0 $.\\
{\it (ii)} Under the assumptions  (\ref{cidn}), (A1)-(A3),(B1)-(B3), (C1)-(C5), we have 
\begin{equation}
\label{eq3.7}
n\|\hat \theta_{2n}-\theta^0_2\|=O_{\eP}(1), \qquad n^{1/2}\|\hat \theta_{1n}-\theta^0_1\|=O_{\eP}(1)
\end{equation}
\end{theorem}

\noindent {\bf Proof of Theorem  \ref{theorem 1}} {\it (i)} To show the strong consistency of the M-estimator, we first prove that  the function $e(\theta):=\eE_{(\varepsilon,X)} \cro{\left| \rho(Y-f_\theta(X))-\rho(Y-f_{\theta^0}(X))\right|}$ is continuous. By the mean value theorem, we have:
\begin{equation}
\label{eq3}
\left|\rho(Y-f_\theta(X))-\rho(Y-f_{\theta^0}(X)) \right| \leq \int^{|\delta_{(\theta,\theta^0)}(X)|}_0 |\psi(\varepsilon-v)|dv
\end{equation}
Then $e(\theta) <\infty$, for all $\theta \in \Omega$, whence the function $e(\theta)$ is well defined. The result of uniform convergence of $e(\theta)$ given by   Lemma \ref{lemma 1} and  $e(\theta^0)=0$ imply that the function $e$ is continuous on $\Omega$. By (\ref{cidn}), we have that $e(\theta)\neq 0$ for all $\theta \neq \theta^0$. Then, we can apply  a similar method to  that in Huber (1967) and  we obtain the strong consistency of the M-estimator.\\
{\it (ii)} Since $\hat \theta_n$ is strongly consistent,  it suffices to suppose $\theta$ in a $\varrho$-neighborhood of $\theta^0$. For a positive constants $b$ and $\varrho$, which will be later determined, let be the sets of parameters:
$$
{\cal  V}_{1b\varrho}:=\acc{\theta \in \Omega_\varrho(\theta^0); \;\; n^{1/2}\|\theta_1-\theta^0_1\|>b}, \quad {\cal V}_{2b\varrho}:=\acc{\theta \in \Omega_\varrho(\theta^0); \;\;n\|\theta_2-\theta^0_2\| >b}
$$
The theorem is proved if we show that: for any $\gamma >0$, $c \in (0,\infty)$ there exist $ b \in (0,\infty)$ and $n_{c \gamma b} \in \eN$ such that:
\begin{equation}
\label{brond}
\eP [\inf_{\theta \in \bar {\cal V}_{jb \varrho}}D_n(\theta_1,\theta_2)>c]>1-\gamma,  \qquad \forall n \geq n_{c \gamma b}, \quad j=1,2
\end{equation}
where $\bar {\cal V}_{jb \varrho}$ is the close of ${\cal V}_{jb \varrho}$. By  relation (\ref{eq3.13}),  we have, for $j=1,2$:
\begin{equation}
\label{eq3.16}
\inf_{\theta \in \bar {\cal V}_{jb\varrho}}D_n(\theta_1,\theta_2) \geq \inf_{\theta \in \bar  {\cal V}_{jb\varrho}}D_n^{(1)} \pth{n^{1/2}(\theta_1-\theta^0_1)}+\inf_{\theta \in \bar  {\cal V}_{jb\varrho}}D_n^{(2)}(\theta_1,\theta_2)
\end{equation}
$\bullet$ The study of $D_n^{(1)}$  is simpler because it involves only the regression parameters:
\begin{equation}
\label{mb}
\inf_{\theta \in \bar  {\cal V}_{2b\varrho}}D_n^{(1)} \pth{n^{1/2}(\theta_1-\theta^0_1)}
 = \min \{\inf_{\| w_1\| \leq b} D_n^{(1)}(w_1) , \inf_{\| w_1\| > b} D_n^{(1)}(w_1) \}
\end{equation}
For $w_1 \in \Gamma^{K+1}$, since $\rho$ is convexe, there are $\tilde b >0$ such that $\inf D_n^{(1)}(w_1)$ is greater than $C\sum^n_{i=1} \rho(\delta_{((\theta^0_1+n^{-1/2}w_1,\theta_2^0),\theta^0)}(X_i))$, for $\|w_1\| >\tilde b$. Assumption (B2) and the convexity of $\rho$ imply that  for all $\gamma,\gamma_1>0$, there are  $b_1 >0$, ($b_1 \geq \tilde b$) and  $ n_\gamma \in \eN$ such that:
\begin{equation}
\label{arond}
\eP \cro{\inf_{\| w_1\| >b_1} D^{(1)}_n(w_1) >\gamma_1} >1-\gamma/2,  \qquad \forall n>n_\gamma
\end{equation}
Using  relation (\ref{arond}) and the approximation of $D_n^{(1)}$ given in  Lemma \ref{lemma 3.4} we obtain that the minimum of (\ref{mb}) is $O_{\eP}(1)$.\\
On the other hand, for $\theta \in \bar  {\cal V}_{1b\varrho}$ we have with an arbitrarily large probability, for $n$ large: $\inf_{\theta \in \bar  {\cal V}_{1b\varrho}}D_n^{(1)} \pth{n^{1/2}(\theta_1-\theta^0_1)}=\inf_{\| w_1\| >b_1} D^{(1)}_n(w_1)$. Relation (\ref{arond}) implies that:\\ $\inf_{\theta \in \bar  {\cal V}_{1b\varrho}}D_n^{(1)} \pth{n^{1/2}(\theta_1-\theta^0_1)}$ is arbitrarily large and positive with a probability close to 1.\\
$\bullet$ We take charge now the study of $D^{(2)}_n$.  For any positive numbers $b$ and $\varrho$, we prove that   $\inf_{\theta \in \bar  {\cal V}_{1b\varrho}}D_n^{(2)}(\theta_1,\theta_2)=O_{\eP}(1)$), using the decomposition:
$$
\inf_{\theta \in \bar  {\cal V}_{1b\varrho}}D_n^{(2)}(\theta_1,\theta_2)=\min\{ \inf_{\theta \in \bar  {\cal V}_{1b\varrho} \cap {\cal V}_{2b\varrho} }D_n^{(2)}(\theta_1,\theta_2),\inf_{\theta \in \bar  {\cal V}_{1b\varrho} \cap {\cal V}^c_{2b\varrho}}D_n^{(2)}(\theta_1,\theta_2) \}
$$
with ${\cal V}^c_{2b\varrho}=\{\theta \in \Omega_\varrho(\theta^0), n\| \theta_2-\theta^0_2\| \leq b   \}$. Taking into account the convexity of $\rho$ and  the approximation on $D_n^{(2)}$ given in   Lemma \ref{lemma 4.1}, we obtain:
$$
\inf_{\theta \in \bar  {\cal V}_{1b\varrho} \cap {\cal V}^c_{2b\varrho}}D_n^{(2)}(\theta_1,\theta_2)=\inf_{\|t\| \leq b, \|w_1\|=b}D^{(2)}_n\pth{\theta^0_1+n^{-1/2}w_1,\theta^0_2+n^{-1}t}=O_{\eP}(1)
$$
Consider $D^{(2)}_n$ for $\theta $ in ${\cal V}_{2 b_2 \varrho}$. By Lemma \ref{l_interm}, for all positive numbers  $ \gamma$ and $c$, it exist $\gamma_2, b_2 \in (0, \infty)$, $\varrho \in (0,1)$, and $n_2 \in \eN$  such that: $\gamma_2b_2 \inf_k \varphi(\tau^0_k) >2c$ and that  relation (\ref{eq3.17}) is true. By (A1), we choose $\varrho \in (0,1)$ sufficiently small such as $\inf_{\tau^0_k \leq x \leq \tau^0_k+\varrho}g(x) \geq g(\tau^0_k)/2$ for all $k=1,...,K$. Then, for $n > b_2/\varrho$ we have: $\inf_{\tau^0_k \leq x \leq \tau^0_k+b_2n^{-1}}g(x) \geq g(\tau^0_k)/2$, for all $k=1,...,K$ and:
\begin{equation}
\label{dn2c}
\eP \cro{\inf_{\theta \in \bar {\cal V}_{2 b_2 \varrho} } D^{(2)}_n(\theta_1,\theta_2)>c} \geq \eP \cro{\inf_{\theta \in \bar {\cal V}_{2 b_2 \varrho}} \frac{D^{(2)}_n(\theta_1,\theta_2)}{nG(|\theta_2-\theta^0_2 |)} > \gamma_0} >1-\frac \gamma 2, \quad \textrm{for } n >n_2
\end{equation}
Hence the second term of the right-hand side  of  inequality (\ref{eq3.16}) is arbitrarily large with arbitrarily large probability for sufficiently large $n$. \\
$\bullet$ In conclusion, we showed that for every set ${\cal V}_{1 b_2 \varrho}$, ${\cal V}_{2 b_2 \varrho}$,  the right-hand side of (\ref{eq3.16}) is the sum of $O_{\eP}(1)$ and of arbitrarily large random variables. This implies  relation (\ref{brond}).
 \hspace*{\fill}$\Diamond$ \\

\subsection{Asymptotic distributions}
We mean now giving the limiting distribution of the M-estimator and an asymptotic approximation for the M-process. \\
Let us consider $t \in \eR^{*K}$ and  $w_1 \in \Gamma^{K+1}$.  For  $D_n$ defined by (\ref{dntheta}) as a process in the standardized parameters, we have  the following decomposition:
\begin{equation}
\label{a}
D_n\pth{\theta^0_1+n^{-1/2}w_1,\theta^0_2+n^{-1}t}=D^{(1)}_n(w_1)+D^{(2)}_n\pth{\theta^0_1+n^{-1/2}w_1,\theta^0_2+n^{-1}t}
\end{equation}
Let us denote $V_0:=\eE_X \cro{\ef_{\theta^0}(X) \ef_{\theta^0}(X)^t}$ the Fisher information matrix corresponding to the random model in $X$. We suppose that the matrix $V_0$ is inversible. \\
The M-process is rescaleted in $D_n$ with regard to the rate of convergence. Let us consider  the random vector:
$$
{\cal Z}_n:=n^{-1/2}\sum^n_{i=1}\ef_{\theta^0}(X_i)\psi(\varepsilon_i)
$$
\noindent Let $D(-\infty,\infty)$ be the set of all cadlag functions on $(-\infty,\infty)$ with the Skorokhod topology.
 
The next theorem gives the joint asymptotic distributions of the M-estimators. In the asymptotic behaviour of regression parameters estimator, the independence of error $\varepsilon$ and of regressor $X$ intervenes in an essential way in variance formula. Also for $\hat \theta_{1n}$, the asymptotic approximation expression is similar to that of the M-estimator in a model without break. On the other hand, the asymptotic distribution of the change-points estimators depends only on the density of $X$ in the true break points and on the difference $\rho(\varepsilon \pm d^0_k)-\rho(\varepsilon)$. 
\begin{theorem}
\label{theorem 3}
Under assumptions  (\ref{cidn}), (A1)-(A3), (B1)-(B3), (C1)-(C5), we have
\begin{equation}
\label{eq4.1}
n^{1/2}(\hat \theta_{1n}-\theta^0_1)=\cro{\el(0)}^{-1} V^{-1}_0{\cal Z}_n+o_{\eP}(1)
\end{equation}
Moreover, $\pth{n^{1/2}(\hat \theta_{1n}-\theta^0_1),n(\hat \theta_{2n}-\theta^0_2)}  \overset{\cal L} {\underset{n \rightarrow \infty}{\longrightarrow}} ({\cal Z}, \Pi_-)$, with\\ ${\cal Z} \sim {\cal N}_{(K+1)d}\pth{0,\eE_\varepsilon\cro{\psi^2(\varepsilon)}\el(0)^{-2} V^{-1}_0}$ a Gaussian random vector  independent of $\Pi_-=\pth{\Pi_{1-}, ...,\Pi_{K-}}$, $\Pi_{k-}=\arg \min_{t_k \in \eR} {\cal P}_k(t_k)$, where:
\begin{equation}
\label{pp}
{\cal P}_k(t_k)={\cal P}_{k1}(t_k)\e1_{t_k \geq 0}+{\cal P}_{k2}(-t_k)\e1_{t_k \leq 0}
\end{equation}
${\cal P}_{k1}$ and ${\cal P}_{k2}$ are two independent compound Poisson processes on $[0,\infty)$ with rate $\varphi(\tau^0_k)$ and ${\cal P}_{k1}(0)={\cal P}_{k2}(0)=0$. The distribution of jumps is given by: $\rho(\varepsilon+d^0_k)-\rho(\varepsilon)$, respectively $\rho(\varepsilon-d^0_k)-\rho(\varepsilon)$.
\end{theorem}

\noindent {\bf Proof of Theorem   \ref{theorem 3}} Using the approximation results obtained in   Lemmas \ref{lemma 3.4}  and  \ref{lemma 4.1} and also  the decomposition: $
M_n\pth{\theta^0_1+n^{-1/2}w_1,\theta^0_2+n^{-1}t}=M_n(\theta^0_1,\theta^0_2)+D_n\pth{\theta^0_1+n^{-1/2}w_1,\theta^0_2+n^{-1}t}$, we obtain an asymptotic approximation for the standardized M-process as the sum of two processes. The first is  quadratic form $Q_n(w_1)$ in the standardized  regression parameters, the second is a empirical process in the standardized change-point parameters:
\begin{equation}
 \label{theorem 4}
M_n\pth{\theta^0_1+n^{-1/2}w_1,\theta^0_2+n^{-1}t}=Q_n(w_1)+D^{(2)}_n(\theta^0_1,\theta^0_2+n^{-1}t)+o_{\eP}(1)
\end{equation}
where:
\begin{equation}
\label{QW}
Q_n(w_1)=M_n(\theta^0_1,\theta^0_2)-n^{-1/2}w_1^t \sum^n_{i=1} \ef_{\theta^0}(X_i)\psi(\varepsilon_i)+\frac{\el(0)}{2} w^t_1 V_0 w_1
\end{equation}
 Let us remind that $\psi$ is the derivative of the function $\rho$.  \\
For $t=(t_1,...,t_K) \in \eR^{*K}$, $w_1 \in \Gamma^{K+1}$, by  relation  (\ref{theorem 4}) we have that the minimum of $M_n\pth{\theta^0_1+n^{-1/2}w_1,r+n^{-1}t}$ with respect to $(w_1,t)$ is equivalent with the minimum with respect to $w_1$ of $Q_n(w_1)$ and with respect to $t$ of $D^{(2)}_n\pth{\theta^0_1,\theta^0_2+n^{-1}t}$. Then relation (\ref{eq4.1}) results from (\ref{QW}). Relation (\ref{eq4.1}) implies that the 
study of the distribution limit of $n^{1/2}(\hat \theta_{1n}-\theta^0_1)$ amount to study the law limit of ${\cal Z}_n$.  But taking into account (C1), by a Central Limit Theorem, ${\cal Z}_n$ converges in distribution to the gaussian distribution: ${\cal N}_{(K+1)d}\pth{0,V_0\eE_\varepsilon[\psi^2(\varepsilon)]}$.\\
In view of Theorem \ref{theorem 1} ({\it ii}) for the change-point estimator, we have:
$$
n(\hat \theta_{2n}-\theta^0_2)=\arg \min_{t  \in \eR^K}D^{(2)}_n\pth{\theta^0_1,\theta^0_2+n^{-1}t}+o_{\eP}(1)
$$
\noindent  For study jointly the distribution of ${\cal Z}_n$ and of $D^{(2)}_n$  we apply  Theorem 4.2 of Koul et al. (2003) for $f_n(X,\varepsilon):=\ef_{\theta^0_1}(X)\psi(\varepsilon)$ and $h_n(X,\varepsilon):=\rho(\varepsilon+d_{(\alpha^0_k,\alpha^0_{k-1})}(X))-\rho(\varepsilon)$. Note that:\\ $
D^{(2)}_n\pth{\theta^0_1,\theta^0_2+n^{-1}t}=\sum^n_{i=1} \sum^K_{k=1}h_n(X_i,\varepsilon_i) \e1_{\min\{\tau^0_k, \tau^0_k+t_k/n\} <X_i \leq \max\{\tau^0_k, \tau^0_k+t_k/n\}}$. On the other hand, for $\xi_n(x,z):=\eE_\varepsilon\cro{\exp\pth{in^{-1/2}z^tf_n(X,\varepsilon)}| X=x}$ we have:
$$
\left| n \pth{1-\xi_n(x,z)} \right| \leq \frac{\eE_\varepsilon\cro{\psi^2(\varepsilon)}}{2} \cro{z^t\ef_{\theta^0_1}(x)}^2 \leq C {\eE_\varepsilon\cro{\psi^2(\varepsilon)}} \|z\|^2 \sup_\alpha \left\|  \frac{\partial h_\alpha(x)}{\partial \alpha}\right\|^2
$$
By  assumptions (\ref{eth}) and (A1), we obtain that $n \pth{1-\xi_n(x,z)}$ is uniformly integrable with respect to $dH(x)$, where $H$ is the distribution function of $X$. Thus $n \pth{1-\xi_n(x,z)}  \rightarrow \eE_{\varepsilon}\cro{\psi^2(\varepsilon)} z^t \Lambda(x)z$,  with: $\Lambda(x):=\ef_{\theta^0_1}(x)\ef^t_{\theta^0_1}(x)$ and $\Lambda:=V_0=\eE_X[\Lambda(X)]$. Whence:
$$
\pth{{\cal Z}_n,D^{(2)}_n({\theta^0_1,\theta^0_2+n^{-1}t})}\overset{{\cal L}} {\underset{n \rightarrow \infty}{\longrightarrow}} \pth{{\cal N}_{(K+1)d}\pth{0,V_0\eE_\varepsilon[\psi^2(\varepsilon)]}, {\cal P}(t)}
$$
in $\eR^{(K+1)d}\times D(-\infty,\infty)^K$ with  ${\cal P}(t):=\sum^K_{k=1}{\cal P}_k(t_k)$. The random vector ${\cal N}_{(K+1)d}\pth{0,V_0\eE_\varepsilon[\psi^2(\varepsilon)]}$ is independent of ${\cal P}_k$, $k=1,...,K$.\\
 We prove now that $n(\hat \theta_{2n}-\theta^0_2)$ converges weakly to the smallest minimizer $\Pi_-$ of the process ${\cal P}$ and show then that the components of this vector coincide with the minimizer of ${\cal P}_k(t_k)$, with the probability 1. Seen the Skorokhod space definition, $D(-\infty, \infty)$, we consider that change-points vary in a compact of $\eR^K$. \\
We consider the M-estimator of the change-points: $\hat \theta_{2n}^b:=\arg \min_{t \in [-b,b]^K}M_n\pth{\hat \theta_{1n}(t),t} $ and the minimizer of ${\cal P}(t)$: $\Pi^b_-:=\arg \min_{t \in [-b,b]^K}{\cal P}(t)$, for a fixed $b> 0$. By  Theorem \ref{theorem 1}, there is a real number $b < \infty$ such that $\hat \theta_{2n} -\hat \theta_{2n}^b \rightarrow 0$ a.s. for $n \rightarrow \infty$. More, it also exists a real $b < \infty$ such that $\Pi_-=\Pi_{-}^b$ with a probability arbitrarily large.\\ Then, we shall first prove that for all $b >0$: 
\begin{equation}
\label{conv}
n(\hat \theta^b_{2n}-\theta^0_2) \overset{{\cal L}} {\underset{n \rightarrow \infty}{\longrightarrow}} \Pi_{-}^b
\end{equation} 
For $t \in [-b,b]^K$, $\tilde b=(b,...,b)$ a K-vector, we consider the random  process ${\cal P}^b(t):={\cal P}(t) \e1_{|t| < \tilde b}$ and: $
\hat M^b_n(t):=\cro{M_n\pth{\hat \theta_{1n}(\theta^0_2+n^{-1}t),\theta^0_2+n^{-1}t}-M_n\pth{\hat \theta_{1n}(\theta^0_2),\theta^0_2}}\e1_{|t| <\tilde b}$. Let also, for $ v \in \eR$, the random process:
$$
H^k_n(v)=\sum^n_{i=1}\cro{\rho(\varepsilon_i+sign(v)d^0_k)-\rho(\varepsilon_i)} \e1_{\min \pth{\tau^0_k,\tau^0_k+n^{-1}v} <X_i  \leq \max \pth{\tau^0_k,\tau^0_k+n^{-1}v}}
$$
and theirs sum: $H_n(t)= \sum^K_{k=1} H^k_n(t_k)$. So by (C3), $\eE_{(\varepsilon,X)} \cro{\sup_{\|t\| \leq b} \left| D^{(2)}_n(\theta^0_2+n^{-1}t)-H_n(t)\right|}$ is bounded to upper by 
$$
 n \sum^K_{k=1} \int_{|x-\tau^0_k| \leq n^{-1}b} \varphi(x) \eE_\varepsilon\cro{\left| \rho\pth{\varepsilon+d_{(\alpha^0_k,\alpha^0_{k-1})}(x)}-\rho\pth{\varepsilon+sign(t_k)d^0_k} \right|} dx
$$
$$
=n \sum^K_{k=1} \int_{|x-\tau^0_k| \leq n^{-1}b} \varphi(x) \eE_\varepsilon\cro{| d_{(\alpha^0_k,\alpha^0_{k-1})}(x)-d^0_k| \left|  \psi(\varepsilon+y_x)\right|} dx, \quad \textrm{with }y_x \rightarrow 0, \textrm{ for } x \rightarrow 0
$$
$$
\leq C n \sum^K_{k=1} \int_{|x-\tau^0_k| \leq n^{-1}b} |x-\tau^0_k|\pth{ \sup_\alpha \left\| \frac{\partial h_\alpha(x)}{\partial \alpha} \right\|} \varphi(x) \eE_\varepsilon\cro{\psi(\varepsilon+y_x)}dx
$$
But $\varphi(x)<C$ and $\eE_\varepsilon\cro{\psi(\varepsilon+y_x)} <C$ as a continuously function on a compact. Then, by the Cauchy-Schwarz inequality:
$$
\eE_{(\varepsilon,X)} \cro{\sup_{\|t\| \leq b} \left| D^{(2)}_n(\theta^0_2+n^{-1}t)-H_n(t)\right|} \leq  C n \sum^K_{k=1} \cro{\int_{|x-\tau^0_k| \leq n^{-1}b} (x-\tau^0_k)^2 dx}^{1/2}=o(1)
$$
Hence: $\sup_{\|t\| \leq  b} \left| \hat M_n^b(t)-H_n(t)\right| =o_{\eP}(1)$. 
Let us consider: $\Pi^b_n=\arg \min_{t \in [-b,b]^K}H_n(t)$. By  Lemmas 4.3 and 4.4  of Koul et al. (2003) we obtain: \\$
n(\hat \theta^b_{2n}-\theta^0_2 ) -  \Pi^b_n \overset{{\cal L}} {\underset{n \rightarrow \infty}{\longrightarrow}} 0$,  $ \Pi^b_n \overset{{\cal L}} {\underset{n \rightarrow \infty}{\longrightarrow}} \Pi^b_-$. 
Then  relation (\ref{conv}) follows. Because for two different change-points we have to make of two independent sets of random variables we have that: $\arg \min_{t \in [-b, b]^K}H_n(t))=\sum^K_{k=1}\arg \min_{t_k \in [-b, b]}H_n^k(t_k))$. The last relation, with (\ref{conv}) and  $\Pi^b_n \overset{{\cal L}} {\underset{n \rightarrow \infty}{\longrightarrow}} \Pi^b_-$,  imply that the asymptotic distribution of $n(\hat \theta_{2n}-\theta^0_2)$ is $\Pi_-$.
\hspace*{\fill}$\Diamond$ \\

{\bf Remarks.} {\bf 1.} For $K=1$, we find the results that the empirical processes $D^{(2)}_n\pth{\theta^0_1,\theta^0_2+n^{-1}t}$ converges to a compound Poisson  process. We also find all the asymptotic distributions for particular estimators ML, LS, LAD.  Particularly, the asymptotic variance of the ML estimator of the regression parameters is: $\eE_\varepsilon[(\psi_\varepsilon'(\varepsilon)/\psi_\varepsilon(\varepsilon))^2]V_0 $, with $\psi_\varepsilon$ the density of $\varepsilon$.  \\
{\bf 2.} Consequence of Theorem \ref{theorem 3}, we can find the  confidence interval  or make   hypothesis test for the parameter $\theta$.\\
{\bf 3.} The discontinuity in the change-points of the regression functions influences the rate of convergence of the change-point estimator.  The proved results are differently from those in the continuous or discontinuous in the change-points for non-random design cases. For example, Van der Geer (1988) prove that in the uniform non-random design two-phase, discontinuous, the limiting distribution of the change-point estimator is determined by a Brownian motion with a linear drift. Rukhin and Vajda (1997) for a continuous model prove that the M-estimator of the change-point is asymptotically normal.

\section{Appendix: Lemmas}
To begin, we state a elementary lemma.
\begin{lemma}
\label{lemma i}
For any $k$ random variables $Z_1,...,Z_k$ the following inequalities are valid:
$$
\sum^k_{i=1}\eP \cro{Z_i<0}-(k-1) \leq  \eP\cro{\pth{\sum^k_{i=1} Z_i}<0} \leq \sum^k_{i=1}\eP \cro{Z_i<0}
$$
\end{lemma}
 
 The following lemma of uniform convergence will be useful in the proofs of the main theorems.
\begin{lemma}
\label{lemma 1} Under assumptions  (A1), (B1), (B3) and  (C2),  we have 
$$
\lim_{\eta \searrow 0}\eE_{(\varepsilon,X)} \cro{\sup_{\theta^* \in  \Omega_\eta(\theta)}\left|\rho\pth{Y-f_\theta(X)}-\rho\pth{Y-f_{\theta^*}(X)} \right|}= 0
$$
\end{lemma}

\noindent {\bf Proof of Lemma \ref{lemma 1}} 
We apply a version of  the mean value theorem:
\begin{equation}
\label{rho}
\rho\pth{Y-f_\theta(X)}-\rho\pth{Y-f_{\theta^*}(X)} =\delta_{(\theta,\theta^*)}(X) \int^1_0 \psi\pth{Y-f_\theta(X)+v\delta_{(\theta,\theta^*)}(X)}dv
\end{equation}
We begin by showing that:
\begin{equation}
\label{exd}
\eE_X \cro{\sup_{(\theta,\theta^*) \in  \Omega_\eta(\theta)} \delta^2_{(\theta,\theta^*)}(X)} {\underset{\eta \rightarrow 0}{\longrightarrow}} 0
\end{equation}
Regarding the change-points,  there are two possible cases. \\
{\bf Case 1}. $\tau_k \in \eR$, $\forall k=1,...,K$. We have:
$$
\sup_{\theta^* \in  \Omega_\eta(\theta)}\left| \delta_{(\theta,\theta^*)}(X) \right| \leq C \cro{\eta  \left\| \sup_{\alpha \in \Gamma} \frac{\partial h_\alpha(X)}{\partial \alpha}\right\|+2 \sup_{\alpha \in \Gamma} | h_\alpha(X)| \sum^K_{k=1} \e1_{|X-\tau_k| \leq \eta}}
$$
Furthermore $\eP \cro{|X-\tau_k| \leq \eta} \rightarrow 0$ for $\eta \rightarrow 0$. Then, with  condition (\ref{cidn}) we obtain (\ref{exd}).\\
{\bf Case 2}. $\tau_1 = -\infty $ or $\tau_K=\infty$. Without loss of generality, we consider $\tau_1 = -\infty $. Obviously $\tau^*_1 \geq\tau_1$. We have $|\tau^*_1-\tau_1| \leq \eta$. Then:
$$
\sup_{\theta^* \in \Omega_\eta(\theta)}| \delta_{(\theta,\theta^*)}(X)| \leq C \cro{\eta  \left\| \sup_{\alpha \in \Gamma} \frac{\partial h_\alpha(X)}{\partial \alpha}\right\|+2 \sup_{\alpha \in \Gamma} | h_\alpha(X)| \pth{\e1_{X <\tau^*_1}+\sum^K_{k=2} \e1_{|X-\tau_k| \leq \eta}}}
$$
But $\eP\cro{X < \tau^*_1} \rightarrow 0$ for $\eta \rightarrow 0$. Using assumption (B3) and the Cauchy-Schwarz inequality, we obtain  relation (\ref{exd}).\\
On the other hand, using the inequality: $\forall x \in \eR$, $
|\psi(x+\varepsilon)| \leq |\psi(\varepsilon+|x|)|+ |\psi(\varepsilon-|x|)|$, we have, by  the Cauchy-Schwarz inequality:
$$
\eE_{(\varepsilon,X)}[\sup_{\theta^* \in \bar \Omega_\eta(\theta)} \left| \delta_{(\theta,\theta^*)}(X)\right| \int^1_0 |\psi(Y-f_\theta(X)+v\delta_{(\theta,\theta^*)}(X))|dv]
$$
$$
\leq C \eE_{(\varepsilon,X)}^{1/2}[\sup_{\theta^* \in \bar \Omega_\eta(\theta)}\delta^2_{(\theta,\theta^*)}(X)]\eE_{(\varepsilon,X)}^{1/2}[\sup_{\theta^* \in \bar \Omega_\eta(\theta)} \psi^2 (\varepsilon+1/2|\delta_{(\theta,\theta^*)}(X)|+ |\delta_{(\theta,\theta^0)}(X)|)]
$$
$$
+C \eE_{(\varepsilon,X)}^{1/2}[\sup_{\theta^* \in \bar \Omega_\eta(\theta)}\delta^2_{(\theta,\theta^*)}(X)]\eE_{(\varepsilon,X)}^{1/2}[\sup_{\theta^* \in \bar \Omega_\eta(\theta)} \psi^2 (\varepsilon-1/2|\delta_{(\theta,\theta^*)}(X)|- |\delta_{(\theta,\theta^0)}(X)|)]
$$
The conclusion results from  relations  (\ref{rho}), (\ref{exd}) and from  assumption (C2).
\hspace*{\fill}$\Diamond$ \\

For $x,z \in \eR$, $\tau \in \eR^K$,  for each $k=1,...,K$,  let be  function: $\nu_k(x,z):=\rho({z+sgn(\tau_k-\tau^0_k) d_{(\alpha^0_k,\alpha^0_{k-1})}(x)})-\rho(z)$.  Let be  the function: $
p_k(x):=\eE_\varepsilon [\nu_k(x,\varepsilon)]$. For each break point $\tau^0_k$ we count the number of $X_i$ which fall into the interval $(\tau^0_k,\tau^0_k+|u_k|)$, with $u=(u_1,...,u_K) \in \eR^K $. Let us consider  the functions $G_k, G_{k,n}: \eR^* \rightarrow (0,1]$, where $\eR^*=\eR \setminus \{0\} $:
$$ G_{k,n}(u_k):=n^{-1}\sum^n_{i=1}\e1_{\min(\tau^0_k,\tau^0_k+u_k) < X_i \leq \max(\tau^0_k,\tau^0_k+u_k)}
$$
and its expectation: $G_k(u_k):=\eE_X[\e1_{\min(\tau^0_k,\tau^0_k+u_k) < X_i \leq \max(\tau^0_k,\tau^0_k+u_k)}]$. For all $K$ change-points we define the functions $G, G_{n}: {\eR^*}^K \rightarrow \eR_+$,
$$
G(u):=\sum^K_{k=1}G_k(u_k), \qquad G_n(u):=\sum^K_{k=1}G_{k,n}(u_k)
$$

We present a lemma that states an important property for this functions.

\begin{lemma}
\label{lemma 2}
Under (A1), for each $\gamma >0$, $\eta>0$, there exists a constant $0 <B<\infty$, such that for all $ b \in (0,1)$, and $ n \geq [B/b]+1$, 
$$
\eP [\sup_{B/n<\|u\| \leq b} |\frac{G_n(u)}{G(u)}-1 | < \eta]>1-\gamma
$$
$$
\eP [\sup_{B/n<\|u\| \leq b} |\frac{  \overset{K}{\underset{k=1}{\sum}} Z^{(1)}_{k,n}(u_k)}{G(u)}-1 | < \eta]>1-\gamma
$$
where $Z^{(1)}_{k,n}: \eR^* \rightarrow \eR$, $k=1,...,K$ is defined by:\\ $Z^{(1)}_{k,n}(u_k):=n^{-1}\sum^n_{i=1} \cro{p_k(X_i)-\nu_k(X_i,\varepsilon_i)}\e1_{\min(\tau^0_k,\tau^0_k+u_k) < X_i \leq \max(\tau^0_k,\tau^0_k+u_k)}
$
\end{lemma}

\noindent {\bf Proof of Lemma \ref{lemma 2}} The proof is similar to that of  Lemma 5.1 of Ciuperca and Dapzol  (2008)  using the results for a single change-point (see Lemma 3.2 of  Koul and Qian (2002)).
\hspace*{\fill}$\Diamond$ \\

Let us now given an approximation of the M-process in $\theta^0$ in the direction of the parameters of regression.
\begin{lemma}
\label{lemma 3.4} Under the assumptions  (A1)-(A3), (B1), (B3), (C1), (C4) and (C5), for each  $ b \in (0,\infty)$, we have:
\begin{equation}
\label{supd1}
\sup_{\|w_1 \|<b} |D^{(1)}_n(w_1)+n^{-1/2}w^t_1\sum^n_{i=1}[\ef_{\theta^0}(X_i)\psi(\varepsilon_i)]-\frac{\el(0)}{2} w^t_1 V_0w_1 |=o_{\eP}(1)
\end{equation}

\end{lemma}

\noindent {\bf Proof of Lemma \ref{lemma 3.4}} Using (C1), (C4) and (C5), we have that: $D^{(1)}_n(w_1)$ is equal to: 
$$\sum^n_{i=1} \cro{f_{(\theta^0_1,\theta^0_2)}(X_i)-f_{(\theta^0_1+n^{-1/2}w_1,\theta^0_2)}(X_i) }\psi(\varepsilon_i)$$
$$-\frac 12 \sum^n_{i=1} \cro{f_{(\theta^0_1,\theta^0_2)}(X_i)-f_{(\theta^0_1+n^{-1/2}w_1,\theta^0_2)}(X_i) }^2\psi'(\varepsilon_i)(1+o_{\eP}(1))$$
Since $\lambda(0)=0$ and by the assumption (\ref{eth}) we obtain:
$$
D^{(1)}_n(w_1)=\acc{-n^{-1/2}w^t_1 \sum^n_{i=1} \ef_{\theta^0}(X_i)\psi(\varepsilon_i)-\frac 12 n^{-1}w^t_1 \sum^n_{i=1} \ef_{\theta^0}(X_i)\ef^t_{\theta^0}(X_i) w_1 \psi'(\varepsilon_i)}(1+o_{\eP}(1))
$$
with $o_{\eP}(1))$ uniformly in  $w_1$ and $n$. Thus, using (A1)-(A3) and (C4), by the strong law of large numbers for $\overset{n} {\underset{i=1}{\sum}}\ef_{\theta^0}(X_i)\ef^t_{\theta^0}(X_i)\psi'(\varepsilon_i)$ and by the assumption (\ref{eth}) for $h_\alpha$, we get:\\ $
D^{(1)}_n(w_1)=\{-n^{-1/2}w^t_1 \ef_{\theta^0}(X_i)\psi(\varepsilon_i)-\frac{\el(0)}{2}w_1V_0 w^t_1\}(1+o_{\eP}(1))$. Thus the proof is complete.
\hspace*{\fill}$\Diamond$ \\

In the following lemma, the set ${\cal V}_{2 b_2 \varrho}$ is defined in the proof of Theorem \ref{theorem 1}. The proof of the similar result for two-phase linear model of Koul et al. (2003) is facilitated by the existence of a single change-point and especially by the linearity in $x$ of $h_\alpha(x)$.
\begin{lemma}
\label{l_interm}
Under the assumptions  (\ref{cidn}), (A1)-(A3), (B1), (B3),  (C1)-(C5), for all positive numbers  $ \gamma$ and $c$, it exist $\gamma_2, b_2 \in (0, \infty)$, $\varrho \in (0,1)$, and $n_2 \in \eN$  such that: $\gamma_2b_2 \inf_k \varphi(\tau^0_k) >2c$ and that:
\begin{equation}
\label{eq3.17}
\eP\cro{\inf_{\theta \in \bar {\cal V}_{2 b_2 \varrho}} \frac{D^{(2)}_n(\theta_1,\theta_2)}{nG(| \theta_2-\theta^0_2|)}>\gamma_2}>1-\gamma/2, \qquad \forall n >n_2
\end{equation} 

\end{lemma}
\noindent {\bf Proof of Lemma \ref{l_interm}}   Let us introduce some notations for ease of exposition.  For each change-point $\tau^0_k$, consider the processes:\\
$S^{(1)}_{k,n}(\theta_1,u_k):=n^{-1} \sum^n_{i=1}\cro{\rho(\varepsilon_i+d_{(\alpha^0_{k},\alpha_{k-1})}(X_i)) - \rho(\varepsilon_i+d_{(\alpha^0_{k},\alpha^0_{k-1})}(X_i)) } \e1_{\min(\tau^0_k,\tau^0_k+u_k)<X_i \leq \max(\tau^0_k,\tau^0_k+u_k)}$\\
$S^{(2)}_{k,n}(\theta_1,u_k):=n^{-1} \sum^n_{i=1}[\rho(\varepsilon_i)-\rho(\varepsilon_i+d_{(\alpha^0_k,\alpha_k)}(X_i))]\e1_{\min(\tau^0_k,\tau^0_k+u_k)<X_i \leq \max(\tau^0_k,\tau^0_k+u_k)}$\\
and the functions: $Z^{(2)}_{k,n}: \eR^* \rightarrow \eR$, $k=1,...,K$:\\
$Z^{(2)}_{k,n}(u_k):=n^{-1} \sum^n_{i=1}\cro{p_k(X_i)-\tau^0_k(r_k)}\e1_{\min(\tau^0_k,\tau^0_k+u_k)<X_i \leq \max(\tau^0_k,\tau^0_k+u_k)}
$\\
Let us consider  $\theta_2=\theta^0_2+u$ with $u=(u_1,...,u_K)$. Given these notations, we see that $n^{-1}D^{(2)}_n(\theta_1,\theta_2)$ can be written as:
 \begin{equation}
\label{eq3.18}
\begin{array}{c}
n^{-1}D^{(2)}_n(\theta_1,\theta_2)=\sum^K_{k=1}p_k(\tau^0_k)G_k(u_k)+\sum^K_{k=1}p_k(\tau^0_k)\cro{G_{k,n}(u_k)-G_k(u_k)}\\
+\sum^K_{k=1}\cro{Z^{(1)}_{k,n}(u_k)+Z^{(2)}_{k,n}(u_k)+S^{(1)}_{k,n}(\theta_1,u_k)+S^{(2)}_{k,n}(\theta_1,u_k)}
\end{array}
\end{equation}
We shall prove that the supremum on the set $\bar {\cal V}_{2 b_2 \varrho}$ of all terms on the right-hand side of (\ref{eq3.18}), except the first, divided by $G(|\theta_2-\theta^0_2|)$ is $o_{\eP}(1)$. On the other hand, we prove that the first term is strictly positive with the probability 1. Remember that $\lambda(y)=\eE_\varepsilon[\rho'(\varepsilon+y)]$.  By Fubini's lemma, by (A2), (C1) and (\ref{cidn}),  we obtain that: 
$$p_k(\tau^0_k)=\eE_\varepsilon\cro{\rho(\varepsilon+d^0_k)-\rho(\varepsilon)}=\int^{\max(0,d_k^0)}_{\min(0, d_k^0)} \lambda(z) dz$$
Since the function $\lambda$ is strictly increasing and $\lambda(0)=0$, we obtain: $
p_k(\tau^0_k)>0$, for each $k=1,...,K$. For all $u_k \leq \varrho$, $k=1,...,K$ we have: $| Z^{(2)}_{k,n}(u_k) | \leq \sup_{0 \leq \upsilon \leq \varrho} \left| p_k(\tau^0_k+\upsilon)-p_k(\tau^0_k) \right| G_{k,n}(u_k)$. 
By assumption (C1), for $n \rightarrow \infty$ and $\varrho \rightarrow 0$, we have:
 $\sup_{B/n <\upsilon \leq \varrho} \left| p_k(\tau^0_k+\upsilon)- p_k(\tau^0_k)\right|=o_{\eP}(1)$ and with  Lemma \ref{lemma 2}  for all $\varrho >0$, there is a $B_1>0$ such that for 
$u_k \in (B/n,\varrho)$,   $k=1,...,K$, and for $n > B_1/\varrho$,  we have: $
0 \leq {G_{k,n}(u_k)}/{G(u)} \leq {G_n(u)}/{G(u)}=1+o_{\eP}(1)$. Hence:
\begin{equation}
\label{Z2K}
\sup_{B_1/n <\|u \| \leq \varrho} \frac{\overset{K} {\underset{k=1}{\sum}} \left| Z^{(2)}_{k,n}(u_k)\right|}{G(u)}=o_{\eP}(1) \qquad \textrm{for } n \rightarrow \infty,\;\;  \varrho \searrow 0
\end{equation}
We have a similar relation for $Z^{(1)}_{k,n}$, for a $B_2>0$ and $n>B_2/\varrho$.\\
For $S^{(1)}_{k,n}$, by  Theorem \ref{theorem 1}(i), for all  $x \in \eR$ and for all  $\theta \in \bar {\cal V}_{2b \varrho}$ we have:  $
| d_{(\alpha_{k-1},\alpha_{k-1}^0)}(x) | <C \varrho$, whence:
$$
| S^{(1)}_{k,n}(\theta_1,u_k)| \leq n^{-1} \sum^n_{i=1}\int^{C\varrho}_{-C\varrho} | \psi (\varepsilon_i+d_{(\alpha^0_k,\alpha^0_{k-1})}(X_i)+v)| dv \e1_{\min(\tau^0_k,\tau^0_k+u_k) <X_i \leq \max(\tau^0_k,\tau^0_k+u_k)}
$$
Applying  Lemma 3.2 of Koul et al. (2003) for: $
J(x,z)=\int^{C\varrho}_{-C\varrho} | \psi\pth{z+d_{(\alpha^0_k,\alpha^0_{k-1})}(x)+v} |dv$ we obtain that there exists a $B_3 \in (0,\infty)$ such that for $n \rightarrow \infty$, $\varrho \searrow 0$:
$$
\sup_{\theta \in \bar {\cal V}_{2B_3\varrho}} \frac{| S^{(1)}_{k,n}(\theta_1,u_k)| }{G_k(u_k)}=o_{\eP}(1)
$$
Then, for all $\gamma >0$, $ \eta > 0$:
$$
\eP [\sup_{\theta \in \bar {\cal V}_{2B_3\varrho}} \frac{\overset{K} {\underset{k=1}{\sum}} | S^{(1)}_{k,n}(\theta_1,u_k)| }{G(u)} <\eta]  \geq \eP [\sup_{\theta \in \bar {\cal V}_{2B_3\varrho}} \sum^K_{k=1} \frac{| S^{(1)}_{k,n}(\theta_1,u_k)|}{G_k(u_k)} < \eta] \geq 1-\gamma
$$
A similar relation holds for $\overset{K} {\underset{k=1}{\sum}}S^{(2)}_{k,n}(\theta_1,u_k)$ for a $B_4>0$. Then, for $ l=1,2$:
\begin{equation}
\label{S1K}
\sup_{\theta \in \bar {\cal V}_{2B_{l+2}\varrho}} \frac{\overset{K} {\underset{k=1}{\sum}}S^{(l)}_{k,n}(\theta_1,u_k)}{G(u)}=o_{\eP}(1), \quad \textrm{for }  \gamma >0, \varrho \searrow 0
\end{equation}
On the other hand, for each $\eta >0$:
$$
\eP [\frac{| \overset{K} {\underset{k=1}{\sum}} p_k(\tau^0_k) [G_{k,n}(u_k)-G_k(u_k)]|}{G(u)} <\eta] \geq \eP [\sum^K_{k=1}\frac{| G_{k,n}(u_k)-G_k(u_k)|}{G_k(u_k)} <\frac{\eta}{\max_k p_k(\tau^0_k)}  ]
$$
$$
 \geq \sum^K_{k=1}\eP [ \frac{| G_{k,n}(u_k)-G_k(u_k)|}{G_k(u_k)} <\frac{\eta}{\max_k p_k(\tau^0_k)}  ] -(K-1)
$$
the last inequality is obtained by  Lemma \ref{lemma i}. By Lemma \ref{lemma 2}, for each $\eta>0$ and $\tilde \gamma >0$ there exists a $B_5>0$ such that the probability  which intervenes in the last inequality is bigger than $1- \tilde \gamma$. Choose $\tilde \gamma=\gamma /(4 K)$ and $\eta <\max_k p_k(\tau^0_k) [4+1/(K \max_k p_k(\tau^0_k))]^{-1} $. 
We  obtain  inequality (\ref{eq3.17}) for $\gamma_2=[\max_k p_k(\tau^0_k)-\eta (4+1/(K \max_k p_k(\tau^0_k)))]/2$, $b_2=\max\{B_1,...,B_5 \}$ and $n_2=b_2/\varrho+1$. \hspace*{\fill}$\Diamond$ \\

 Following  result gives the  behaviour of $D^{(2)}_n$ in a $n^{-1/2}$-neighborhood of $\theta^0_1$.
\begin{lemma}
\label{lemma 4.1} Under  assumptions  (\ref{cidn}), (A1)-(A3), (B1), (B3), (C2),  if we define\\ $A_n(w_1,t):=D^{(2)}_n\pth{\theta^0_1+n^{-1/2}w_1,\theta^0_2+n^{-1}t}- D^{(2)}_n\pth{\theta^0_1,\theta^0_2+n^{-1}t}$, we have, for every $b \in (0,\infty)$:
\begin{equation}
\label{supd2}
\sup_{\|( w_1,t)\| \leq b} | A_n(w_1,t)|=o_{\eP}(1)
\end{equation}
\end{lemma}
\noindent {\bf Proof of Lemma \ref{lemma 4.1}}  Without loss of generality, we consider the vector $t=(t_1,...,t_K)\in \eR^{*^K}_+$.  The general case is obtained by very similar arguments. Let us note $w_1=\pth{w_{1,0}, w_{1,1},...,w_{1,K}} \in \Gamma^{K+1}$. First, observe that we can write 
$$
A_n(w_1,t)=\sum^n_{i=1} \sum^{K+1}_{k=1}\int^{d_{(\alpha^0_{k-1},\alpha^0_{k-1}+n^{-1/2}w_{1,k-1} )}(X_i)}_0 \psi(\varepsilon_i +h_{\alpha^0_k}(X_i)+v)dv \e1_{\tau^0_k<X_i \leq \tau^0_k+n^{-1}t_k}
$$
$$
-\sum^n_{i=1} \sum^{K+1}_{k=1}\int^{d_{(\alpha^0_{k-1},\alpha^0_k+n^{-1/2}w_{1,k} )}(X_i)}_0\psi(\varepsilon_i+v)dv \e1_{\tau^0_k<X_i \leq \tau^0_k+n^{-1}t_k}
$$
Let us denote:
\begin{equation}
\label{BAB}
B(w_1,t):=A_n(w_1,t)+B_1(w_1,t)-B_2(w_1,t)
\end{equation}
with $$B_1(w_1,t):=\sum^n_{i=1} \sum^K_{k=1}d_{(\alpha^0_{k-1}+n^{-1/2}w_{1,k-1},\alpha^0_{k-1})}(X_i)\psi(\varepsilon_i)\e1_{\tau^0_k<X_i \leq \tau^0_k+n^{-1}t_k}$$ and
$$B_2(w_1,t):=\sum^n_{i=1} \sum^{K+1}_{k=1}d_{(\alpha^0_k+n^{-1/2}w_{1,k},\alpha^0_k)}(X_i)\psi(\varepsilon_i)\e1_{\tau^0_k<X_i \leq \tau^0_k+n^{-1}t_k} $$
By assumption (\ref{eth}) for  $ \partial h_\alpha/ \partial \alpha$ we have that: $|d_{(\alpha^0_{k-1}+n^{-1/2}w_{1,k-1},\alpha^0_{k-1})}(X_i)| \leq n^{-1/2}U_i$ with $U_i$ a random variable such as $\eE_X(U_i) < \infty$  and there exists  a real $c_{1,b}>0$ such that $\eP\cro{U_i<c_{1,b}}>1-\gamma$, for all $\gamma>0$. \\
Since  $|h_{\alpha^0_{k-1}}(X_i)| \e1_{\tau^0_k<X_i \leq \tau^0_k+n^{-1}t_k}<C$ with the probability 1, we have that $\eE_{(\varepsilon,X)}[|B_1(w_1,t)|]$ is bounded to upper by :
$$
n \int^{n^{-1/2}c_{1,b}}_0 \acc{\eE_\varepsilon[| \psi(\varepsilon+C+v)-\psi(\varepsilon)|]
+\eE_\varepsilon[| \psi(\varepsilon-C+v)-\psi(\varepsilon)|]}dv 
{\sum^K_{k=1}\eE_X\cro{\e1_{\tau^0_k <X_i \leq \tau^0_k+n^{-1}t_k}}}
$$
Since $\varepsilon$ and $X$ are independently and by the relation  (\ref{n2}): \\
$\eE_{(\varepsilon,X)}[|B_1(w_1,t)|]$  $ \leq  o(1)  n \sum^K_{k=1}G_k(n^{-1}t_k)=o(1)$. Alike:
$$
\eE_{(\varepsilon,X)}[|B_2(w_1,t)|] \leq n \int^{n^{-1/2}c_{1,b}}_0 \eE_\varepsilon\cro{| \psi(\varepsilon+v)-\psi(\varepsilon)|}dv \sum^K_{k=1}G_k(n^{-1}t_k)=o(1)
$$
Therefore, we have  
\begin{equation}
\label{eo1}
\eE_{(\varepsilon,X)}[|B_1(w_1,t)|]=o(1), \qquad \eE_{(\varepsilon,X)}[B_2(w_1,t)]=o(1)
\end{equation}
Let us note the random process: $D_1(w_1,t):=B(w_1,t)-A_n(w_1,t)=B_1(w_1,t)-B_2(w_1,t)$. Then, we can write:
$$
D_1(w_1,t)=\sum^n_{i=1}\sum^{K+1}_{k=1}\{ d_{(\alpha^0_{k-1}+n^{-1/2}w_{1,k-1},\alpha^0_{k-1})}(X_i)-d_{(\alpha^0_k+n^{-1/2}w_{1,k},\alpha^0_k)}(X_i) \}\psi(\varepsilon_i)\e1_{\tau^0_k <X_i \leq \tau^0_k+n^{-1}t_k}
$$
Using assumption (\ref{eth}) for  $ \partial^2 h_\alpha/ \partial \alpha^2$ we obtain:
\begin{equation}
\label{eo2}
D_1(w_1,t)=n^{-1/2}\sum^n_{i=1}\sum^{K+1}_{k=1}[w_{1,k-1}\frac{\partial h_{\alpha^0_{k-1}}(X_i)}{\partial \alpha}- w_{1,k}\frac{\partial h_{\alpha^0_{k}}(X_i)}{\partial \alpha}]\psi(\varepsilon_i)\e1_{\tau^0_k <X_i \leq \tau^0_k+n^{-1}t_k}+o_{\eP}(1)
\end{equation}
Finally, by  assumption (C2) for $c=0$ and by (\ref{n2}), we obtain that $\forall b \in (0,\infty), \forall w_1 \in \eR^{K+1}$:
\begin{equation}
\label{eo3}
\sup_{0 \leq \|t\| \leq b} \|\sum^n_{i=1}\sum^{K+1}_{k=1}[w_{1,k-1}\frac{\partial h_{\alpha^0_{k-1}}(X_i)}{\partial \alpha}- w_{1,k}\frac{\partial h_{\alpha^0_{k}}(X_i)}{\partial \alpha}]\psi(\varepsilon_i)\e1_{\tau^0_k <X_i \leq \tau^0_k+n^{-1}t_k} \| =O_{\eP}(1)  
\end{equation}
The conclusion follows from the relations (\ref{BAB}), (\ref{eo1}), (\ref{eo2}) and (\ref{eo3}) .
\hspace*{\fill}$\Diamond$ \\




\end{document}